\documentclass[12pt]{article}

\usepackage{amsfonts}
\usepackage{xcolor}
\usepackage[english]{babel}
\def\findemo{\hfill \rule{6pt}{6pt}}
\def\R{\mathbb{ R}}

\def\N{\mathbb{ N}}
\def\<{\langle}
\def\>{\rangle}

\begin{document}

\baselineskip=14pt
\frenchspacing

\newtheorem{thm}{Theorem}[section]
\newtheorem{theo}[thm]{Theorem}
\newtheorem{prop}[thm]{Proposition}
\newtheorem{coro}[thm]{Corollary}
\newtheorem{lema}[thm]{Lemma}
\newtheorem{defi}[thm]{Definition}
\newtheorem{ejem}[thm]{Example}
\newtheorem{rema}[thm]{Remark}
\newtheorem{fact}[thm]{Fact}
\newtheorem{open}[thm]{PROBLEM}

\title{Generalized metric properties of spheres\linebreak and renorming of Banach spaces}

\author{S. Ferrari, J. Orihuela, M. Raja}

\maketitle

\section{Introduction}

Throughout this paper $X$ will denote a real Banach space and $X^*$ will denote its topological dual. If $F$ is a subset of $X^*$, then $\sigma(X,F)$ denotes the weakest topology on $X$ that makes each member of $F$ continuous or, equivalently, the topology of pointwise convergence on $F$. Analogously, if $E$ is a subset of $X$, then $\sigma(X^*,E)$ is the topology for $X^*$ of pointwise convergence on $E$. We denote with $B_X$ (respectively, $S_X$) the unit ball of $X$ (respectively, the unit sphere of $X$).
Recall that a subset $B$ of $B_{X^*}$ is said to be \emph{norming} if
\[\|x\|_B=\sup_{b^*\in B}|b^*(x)|\]
is a norm on $X$ equivalent to the original norm of $X$. A subspace $F\subseteq X^*$ is norming if $F\cap B_{X^*}$ is norming.\\

In our previous paper \cite{FOR} we have discussed the possibility of metrizing a weak topology on the unit sphere of an equivalent norm of a Banach space. As our result \cite[Theorem 1.4]{FOR}, given a norming subspace $F\subseteq X^*$, the existence of an equivalent norm on $X$ such that the $\sigma(X,F)$ topology is metrizable on its unit sphere is closely related to the existence of a $\sigma(X,F)$-LUR norm, actually a slightly stronger property can be obtained, namely a $\sigma(X,F)$-lsc equivalent norm $|\!|\!|\cdot |\!|\!|$ exists such that: let $x\in X$ and $(x_n)_{n\in\N},(y_n)_{n\in\N}\subseteq X$, we have $\sigma(X,F)$-$\lim_{n}y_n=x$ whenever
\[\lim_{n\rightarrow+\infty}\left (2|\!|\!|x|\!|\!|^2+2|\!|\!|x_n|\!|\!|^2-|\!|\!|x+x_n|\!|\!|^2\right )+\left( 2|\!|\!|x_n|\!|\!|^2+2|\!|\!|y_n|\!|\!|^2-|\!|\!|x_n+y_n|\!|\!|^2\right)=0.\]
On the other hand, the existence of an equivalent $\sigma(X,F)$-LUR norm was characterized by the second named author in terms of countable decompositions of the space $X$, as well as network conditions close to covering properties of generalized metrizable spaces, (see \cite {cas-ori} and \cite{Orihuela}). Moreover, as we shall see here, it is somehow natural to think that the existence of a $\sigma(X,F)$-LUR renorming could be characterized by covering properties of the sphere closed to its metrizability for the $\sigma(X,F)$-topology.\\

In order to state our results we need to recall some basic definitions coming from Banach space geometry and general topology. A norm defined on the space $X$ is said to be \emph{rotund} (or \emph{strictly convex}) if, given $x,y\in X$ satisfying $\|x\|=\|y\|=\|2^{-1}(x+y)\|$ (equivalently, $2\|x\|^2+2\|y\|^2-\|x+y\|^2=0$), we have $x=y$. Geometrically, this means that the unit sphere of $X$ has no non-trivial line segment, or equivalently ${\rm ext}(B_X)=S_X$, meaning that the extreme points of the unit ball are all the points of the unit sphere.
Let $\tau$ be a topology on a normed space $X$. $X$ is said to be $\tau$-locally uniformly rotund ($\tau$-LUR, for short) if given $x\in X$ and $(x_n)_{n\in\N}\subseteq X$, then $x_n$ converges to $x$ in the $\tau$-topology, whenever
\[\lim_{n\rightarrow+\infty}(2\|x\|^2+2\|x_n\|^2-\|x+x_n\|^2)=0.\]
If $\tau$ is the norm topology, we simply say that $X$ is LUR.\\

Let $T$ be a non empty set. We say that a function $\rho: T\times T \rightarrow [0,+\infty)$ is a symmetric on a set $T$ if $\rho(x,y)=\rho(y,x)$ for every $x,y\in T$ and $\rho(x,y)=0$ if, and only if, $x=y$. Balls $B_\rho(x,\varepsilon)$ are defined in the same fashion that for metrics. Because of a (possible) lack of continuity of the function $\rho$, ``open'' balls with respect to a symmetric may not define a topology, but somehow topologies on $T$ may be related to symmetrics. Indeed a topological space $T$ is said to be symmetrizable if there is a symmetric $\rho$ on it such that: $U \subseteq T$ is open if and only if for every $x \in U$ there is $\varepsilon>0$ such that $B_\rho(x,\varepsilon) \subseteq U$. If the balls $B_\rho(x, \varepsilon)$ are a neighbourhood base at $x$ for every $x \in T$ we say that $T$ is semi-metrized by $\rho$ and $T$ is semi-metrizable.
Let us say that a sequence of open coverings $({\mathcal G}_n)_{n \in {\Bbb N}}$ of a topological space $T$ is a development of $X$ if 
\[{\rm St}({\mathcal G}_n,x):=\bigcup\{U\in\mathcal{G}_n: x\in U\}\]
is a neighbourhood base of the topology at every $x \in T$. A regular topological space is said  Moore if it has a development.\\

Our main result links the Moore property, the existence of symmetrics with good properties and the $\tau$-LUR property, all for different equivalent norms.

\begin{theo}\label{moore-slice}
Let $X$ a Banach space and $F \subseteq X^*$ a norming subspace. The following statements are equivalent:
\begin{itemize}
\item[(i)] $X$ admits an equivalent $\sigma(X,F)$-lsc and $\sigma(X,F)$-LUR norm;
\item[(ii)] $X$ admits an equivalent $\sigma(X,F)$-lsc norm such that $(S_X,\sigma(X,F))$ is symmetrized by a symmetric $\rho$ satisfying that every point $x \in S_X$ is contained in a $\sigma(X,F)$-open slice of arbitrarily small $\rho$-diameter;
\item[(iii)] $X$ admits an equivalent $\sigma(X,F)$-lsc norm such that there is a development  $({\mathcal F}_n)_{n \in {\Bbb N}}$ of $(S_X,\sigma(X,F))$ made up of $\sigma(X,F)$-open slices.
\end{itemize}
\end{theo}

Note that under the hypotheses of statement $(ii)$, the sphere $(S_X,\sigma(X,F))$ is actually semi-metrized by $\rho$ because $x$ is in the interior of $B_\rho(x,\varepsilon)$ for every $x \in S_X$ and $\varepsilon>0$. However, the property about small diameters is not guaranteed by the sole symmetrizability or semi-metrizability as the triangle inequality may fails.\\

The techniques used to prove the previous result also sheds light between the covering properties of spheres and the existence of rotund renormings. We say that a topological space $X$ has a ${\mathcal G}_\delta$-diagonal if the diagonal set 
$\Delta:=\{(x,x)\in X\times X: x \in X\}$ is a ${\mathcal G}_\delta$-set in the product space $X \times X$. The following result shows that the ${\mathcal G}_\delta$-diagonal property is a sort of developability.

\begin{theo}[{\cite{Ced61} and \cite[Theorem 2.2]{Gru84}}]
A space $X$ has a ${\mathcal G}_\delta$-diagonal if, and only if, there exists a sequence $({\mathcal G}_n)_{n\in {\Bbb N}}$  of open covers of $X$ such that for each $x, y \in X$ with $x \not =  y$, there exists $n \in {\Bbb N}$ with 
$y \not \in {\rm St}(x,  {\mathcal G}_n)$.
In this case we say that $({\mathcal G}_n)_{n\in {\Bbb N}}$ is a ${\mathcal G}_\delta$-diagonal sequence for $X$
\end{theo}

With this terminology we can state the following result which links rotund renormability with linear topological properties of spheres for equivalent norms.

\begin{theo}\label{diagonal-slice}
Let $X$ be a Banach space and let $F \subseteq X^*$ be a norming subspace. The following statements are equivalent:
\begin{itemize}
\item[(i)] $X$ admits an equivalent $\sigma(X,F)$-lsc rotund norm;
\item[(ii)] $X$ admits an equivalent $\sigma(X,F)$-lsc norm such that there is a symmetric $\rho$ defined on $S_X$ with the property that every point $x \in S_X$ is contained in a $\sigma(X,F)$-open slice of arbitrarily small $\rho$-diameter;
\item[(iii)] $X$ admits an equivalent $\sigma(X,F)$-lsc norm such that $(S_X,\sigma(X,F))$ has a ${\mathcal G}_\delta$-diagonal sequence whose $\sigma(X,F)$-open sets are slices.
\end{itemize}
\end{theo}

Note the difference between the existence of a symmetric in statement $(ii)$ above and being symmetrized by the symmetyric in statement $(ii)$ of Theorem \ref{moore-slice}. Let us recall that rotundity was characterized in \cite{OST} using a linear version of the topological property $(*)$ that we will introduce later (Definition \ref{stardef}).\\

Now we shall state our results for dual Banach spaces equipped with the weak$^*$-topology. In this case, the compactness allows us to improve the results. Let us gather together the most significative characterizatons of $w^*$-LUR renormings from \cite{FOR} with the new ones obtained here.

\begin{theo}\label{moore-dual}
Let $X^*$ be a dual Banach space. The following statements are equivalent:
\begin{itemize}
\item[(i)] $X^*$ admits an equivalent dual $w^*$-LUR norm;
\item[(ii)] $X^*$ admits an equivalent dual norm such that $(S_{X^*},w^*)$ is a Moore space;
\item[(iii)] $X^*$ admits an equivalent dual norm such that $(S_{X^*},w^*)$ is 
symmetrized by a symmetric $\rho$ satisfying that every point $x^* \in S_{X^*}$ has $w^*$-neighborhoods of arbitrarily small $\rho$-diameter;
\item[(iv)] $X^*$ admits an equivalent dual norm such that $(S_{X^*},w^*)$ is metrizable;
\item[(v)] $(B_{X^*},w^*)$ is a descriptive compact space.
\end{itemize}
\end{theo}

The symmetrization hypothesis in statement $(iii)$ is in that context equivalent to semi-metrization, but it cannot be replaced by the mere existence of a symmetric $\rho$ with the small $\rho$-diameter property. We will see that actually such a weaker property is equivalent to the existence of a rotund dual norm, see Theorem \ref{strict-dual} below.
A different question is to know whether the hypothesis on small diameter sets could be removed. We do not know the answer, but our guess is negative. See Corollary \ref{semi-rotund} in relation to this problem.\\

For the last result we need the following definition introduced in~\cite{OST} in relation with rotund renorming:

\begin{defi}\label{stardef}
A topological space $A$ has property $(\ast)$ if there exists a sequence $\{\mathcal{U}_n: n\in{\Bbb N}\}$ of families of open subsets of $A$ such that, given any $x,y\in A$, there exists $n_0\in\N$ such that:
	\begin{itemize}
	\item[a)] $\{x,y\}\cap \bigcup \mathcal{U}_{n_0}\neq\emptyset$, where $\bigcup \mathcal{U}_{n_0}:=\bigcup\{U: U\in\mathcal{U}_{n_0}\}$;
	\item[b)] for every $U\in\mathcal{U}_{n_0}$ the set $\{x,y\}\cap U$ is at most a singleton.
	\end{itemize}
\end{defi}

If $A$ is a subset of a topological vector space and the elements of $\bigcup_{n=1}^\infty\mathcal{U}_n$ can be taken to be slices of $A$, then $A$ is said to have $(\ast)$ with slices. It is shown in~\cite[Theorem 2.7]{OST} that if $X^*$ is a dual Banach space then $(B_{X^*},w^*)$ has $(\ast)$ with slices if and only if $X^*$ admits a dual rotund norm.\\

In this paper we are able to give characterizations of dual rotund renorming  in terms of topological properties of $(S_{X^*},w^*)$ up to equivalent renormings.

\begin{theo}\label{strict-dual}
Let $X^*$ be a dual Banach space. The following statements are equivalent:
\begin{itemize}
\item[(i)] $X^*$ admits an equivalent dual rotund norm;
\item[(ii)] $X^*$ admits an equivalent dual norm such that $(S_{X^*},w^*)$ has a ${\mathcal G}_\delta$-diagonal;
\item[(iii)] $X^*$ admits an equivalent dual norm such that $(S_{X^*},w^*)$ has a symmetric $\rho$ such that every point $x^* \in S_{X^*}$ has $w^*$-neighborhoods of arbitrarily small $\rho$-diameter;
\item[(iv)] $X^*$ admits an equivalent dual norm such that $(S_{X^*},w^*)$ has $(*)$ with a family of sets $({\mathcal U}_n)_{\in {\Bbb N}}$ such that every $\mathcal{U}_n$ is a covering of $S_{X^*}$;
\item[(v)] $(B_{X^*},w^*)$ has $(*)$ with slices.
\end{itemize}
\end{theo}

Note that statement $(v)$, coming from \cite{OST} is not topological, actually, it is linear-topological. However, it is an isomorphic characterization, meanwhile the previous statements depend on the particular dual norm of the space. From the well known fact that semi-metrizable space has ${\mathcal G}_\delta$-diagonal, see \cite[p. 484]{Gru84}, we deduce immediately the following.

\begin{coro}\label{semi-rotund} 
Let $X^*$ be a dual space endowed with a dual norm such that $(S_{X^*},w^*)$ is semi-metrizable. Then $X^*$ admits an equivalent dual rotund norm.
\end{coro}

The rest paper is organized as follows. In the second section we describe a technique of renorming that allow us to do very precise localization of points on a family of slices. That is applied to the proofs of Theorems \ref{moore-slice} and \ref{diagonal-slice}. The last section is devoted to the results in the dual case, Theorems \ref{moore-dual} and \ref{strict-dual}, which requiere a special treatment.

\section{A slice localisation principle}

When dealing with rotund renormings a useful tool is the so called Lancien's midpoint argument which essentially consists to say that two points $x,y$ in a convex $C$ set are ``close'' if the middle point $(x+y)/2$ is inside one ``small'' slice of the set since one of them, either $x$ or $y$, is in the slice too. Clearly this argument works when closeness and smallness are measured in terms of a norm, but for a more specific use we need to localise the points $x$ and $y$ in the same slice. For that aim we need a slice localisation principle in the sense developed by the second named author together with S. Troyansky in \cite[Theorem 3]{OT1}. The proof there depends on Deville's master lemma for the decomposition method (see Lemma 1.1, Capter 7 in \cite{DGZ}). We shall give an alternative proof here, based on more geometrical ideas, when we restrict ourselves to the case where all points involved are assumed to be in the fixed sliced set.
We will use the following notation: given ${\mathcal F}$ a family of sets, then 
$\bigcup {\mathcal F}:=\bigcup_{H \in {\mathcal F}} H$.

\begin{theo}\label{localisation}
Given ${\mathcal S}$ a family of $\sigma(X,F)$-open slices of a set  $A \subseteq X$, there exists a $\sigma(X,F)$-lsc equivalent norm on $X$ such that whenever $x \in  \bigcup {\mathcal S}$ and $(x_n)_{n\in{\Bbb N}} \subseteq A$ satisfy 
$$ \lim_{n\rightarrow+\infty} (2\|x_n\|^2+2\|x\|^2-\|x_n+x\|^2) =0, $$
then there exists $N \in {\Bbb N}$ such that for every $n \geq N$ there is $S_n \in 
{\mathcal S}$ such that $\{x_n,x\} \subseteq S_n$.
\end{theo}

\begin{rema}
The norm given by Theorem \ref{localisation} can be taken arbitrarily close to the original norm of $X$. Indeed, if $|\!|\!| \cdot |\!|\!|$ has the property of the Theorem \ref{localisation}, then any norm defined as $|\!|\!| x |\!|\!|_\varepsilon^2=\|x\|^2+\varepsilon |\!|\!| x |\!|\!|^2$ has the same property for every $\varepsilon>0$ (where $\|\cdot\|$ is the natural norm of $X$).
For that reason, when handling several of such kind of norms we may always assume that have the same bounds of equivalence with respect to the original norm of $X$.
\end{rema}

The proof needs several intermediate renormings, so for the sake of clarity we will split the construction with the help of the following lemma.

\begin{lema}
Let $X$ be a dual space with predual $F$ and $\varepsilon>0$. If ${\mathcal F}$ is a family of $w^*$-open half-spaces such that $\delta B_X \cap \bigcup {\mathcal F} = \emptyset$, for some $\delta>0$, then there exists an equivalent dual norm $|\!|\!| \cdot |\!|\!|_{{\mathcal F},\varepsilon}$ with the following property: given 
$x \in S_X \setminus((B_X \setminus \bigcup {\mathcal F}) +\varepsilon B_X)$
and $(x_n)_{n\in{\Bbb N}} \subseteq X$ with
$$ \lim_{n\rightarrow+\infty} (2|\!|\!|x_n|\!|\!|_{{\mathcal F},\varepsilon}^2+2|\!|\!|x|\!|\!|_{{\mathcal F},\varepsilon}^2-|\!|\!|x_n+x|\!|\!|_{{\mathcal F},\varepsilon}^2 )=0, $$
there is $N \in {\Bbb N}$ such that for all $n \geq N$ there exists $H_n \in 
{\mathcal F}$ with $\{x_n,x\} \subseteq H_n$.
\end{lema}

\noindent
{\bf Proof.} Without loss of generality we can assume that 
\[{\mathcal F} = \{\{f_\gamma>\alpha_\gamma\} : \gamma \in \Gamma \}\]
where $f_\gamma \in S_F$, $\alpha_\gamma\in\R$ and $\Gamma$ is a non empty set. Put $\beta_\gamma=\alpha_\gamma +\varepsilon$ and observe that if $x \not \in (B_X \setminus \bigcup {\mathcal F}) +\varepsilon B_X$, then $ x \in \bigcup_{\gamma \in \Gamma} 
\{f_\gamma>\beta_\gamma\} $.
For every $k \in {\Bbb N}$ we define a new family of half-spaces 
$$ {\mathcal F}_k = \{ \{ f_\gamma > 2^{-k} + (1-2^{-k})\beta_\gamma \} : \gamma \in \Gamma \}. $$
Consider now the $w^*$-closed symmetric convex bodies $C_0=B_X$ and $C_k=B_X\setminus \bigcup {\mathcal F}_k$ for $k \geq 1$. Let $\| \cdot \|_k$ be the Minkowski functional of $C_k$ and
take 
$$ |\!|\!| x |\!|\!|_{{\mathcal F},\varepsilon}^2 = \sum_{k=1}^{+\infty} 2^{-k}\|x\|^2_k .$$
We will check that $|\!|\!| \cdot |\!|\!|_{{\mathcal F},\varepsilon}$ is the desired norm. Assume that 
$(x_n) \subseteq S_X$ satisfy that
$$ \lim_{n\rightarrow+\infty} (2|\!|\!|x_n|\!|\!|_{{\mathcal F},\varepsilon}^2+2|\!|\!|x|\!|\!|_{{\mathcal F},\varepsilon}^2-|\!|\!|x_n+x|\!|\!|_{{\mathcal F},\varepsilon}^2) =0, $$
for some $x \in S_X$. A standard convexity arguments shows that the limit happens in the same way for all the norms $\| \cdot \|_k$. If $x$ satisfy the statement, then there is $k$ such that $x \in C_{k-1} \setminus C_k$. As $x \in C_{k-1}$, we have $\|x\|_{k-1} \leq 1$ and $\|x\|_k >1$. From the convexity we have
$$ \lim_{n\rightarrow+\infty} \|x_n\|_{k-1} = \lim_{n\rightarrow+\infty} \left\| \frac{x_n+x}{2}\right\|_{k-1} = \|x\|_{k-1},$$
$$ \lim_{n\rightarrow+\infty} \|x_n\|_{k} = \lim_{n\rightarrow+\infty} \left\| \frac{x_n+x}{2}\right\|_{k} = \|x\|_{k}.$$
Take $\lambda_n = \|x\|_{k-1}/\|x_n\|_{k-1}$ and $x'_n=\lambda_n x_n$ and observe
$\|\frac{x'_n+x}{2}\|_{k-1}\leq 1$ and $\lim_n \|\frac{x'_n+x}{2}\|_k=\|x\|_k>1$. 
Now fix $N$ such that for $n \geq N$ we have 
$$ \frac{x'_n+x}{2} \in C_{k-1} \setminus C_k $$
and $|\lambda_n-1|<\varepsilon$.
For such $n$ there is a half-space
$$ H_n=\{ f_\gamma > 2^{-k} + (1-2^{-k})\beta_\gamma \} \in {\mathcal F}_{k}$$ 
such that $\frac{x'_n+x}{2} \in H_n$, that is
$f_\gamma(\frac{x'_n+x}{2}) > 2^{-k} + (1-2^{-k})\beta_\gamma$. 
We also have
$f_\gamma(x) \leq 2^{1-k} + (1-2^{1-k})\beta_\gamma$. 
Indeed, otherwise 
$$x \not \in B_X \setminus \{f_\gamma > 2^{1-k} + 
(1-2^{1-k})\beta_\gamma\} \supset C_{k-1}.$$
The inequalities together give us $f_\gamma(x'_n) > \beta_\gamma$. That implies $f_\gamma(x_n) > \alpha_\gamma$. As $f_\gamma(x)>\beta_\gamma>\alpha_\gamma$, we get that $x_n$ and $x$ belong to $\{f_\gamma>\alpha_\gamma\} \in {\mathcal F}$.\findemo\\

\noindent
{\bf Proof of Theorem \ref{localisation}.}
Firstly we may assume that $A$ is bounded. Otherwise it can be decomposed into countable many bounded pieces. The norm will be obtained as a weighted series of squared norms having the property of the statement built for each piece.\\
Now consider the space $F \bigoplus {\Bbb R}$ and note that $X$ embeds isometrically as an affine hyperplane of $F^* \bigoplus {\Bbb R}$, particularly into $F^* \times 
\{1\}$. Clearly, the $\sigma(X,F)$-topology of $X$ is still induced by the $w^*$-topology, and the elements of the family ${\mathcal F}$ may be replaced by $w^*$-open half-spaces in $F^* \bigoplus {\Bbb R}$ which skip 
$B=(1/2)B_{F^* \bigoplus {\Bbb R}}$.\\
There is an equivalent dual norm $p$ on $F^* \bigoplus {\Bbb R}$ whose unit ball contains $B$ and whose unit sphere contains $A$. Indeed, just let $B_p$ the $w^*$-closed convex hull of $A \cup B$.
We have 
$$ A \subseteq \bigcup_{\varepsilon>0} S_p \setminus((B_p \setminus \bigcup {\mathcal F}) +\varepsilon B_p)$$
Let $|\!|\!| \cdot |\!|\!|_k$ the norm given by the lemma for $\varepsilon=1/k$ and observe that the norm
$$ |\!|\!| x |\!|\!|_0^2= \sum_{k=1}^{+\infty} 2^{-k}|\!|\!| x |\!|\!|_k^2 $$
has the desired property, but its restriction to $X$ may not be a norm. Indeed it is just a $\sigma(X,F)$-lsc convex function that we may symmetrize taking 
$$F(x)= |\!|\!| (x,1) |\!|\!|_0^2 + |\!|\!| (-x,1) |\!|\!|_0^2.$$
To get a norm on $X$ we proceed as follows.
Let $m \geq 0$ be the infimum of $F$ and let $\{r_k: k \in{\Bbb N}\}$ be an enumeration of ${\Bbb Q}\cap(m,+\infty)$. Let $\| \cdot \|_k$ the Minkowski functional of the $w^*$-closed convex symmetric body $\{x \in X: F(x) \leq r_k\}$, and finally take
$$ \|x\|^2= \sum_{k=1}^{+\infty} a_k \|x\|_k^2 $$
where the numbers $a_k>0$ are chosen to guarantee the uniform convergence of the series on bounded subsets of $X$. The usual convexity arguments show that this norm has the desired property.\findemo\\

\noindent
{\bf Proof of Theorems \ref{moore-slice} and \ref{diagonal-slice}.} 
Due to the similarity between the statements we will proceed in ``parallel'' with both proofs.
Assume $(i)$ so $X$ is endowed with a $\sigma(X,F)$-lsc rotund norm. Note that
$$ \rho(x,y) = 1 - \left\| \frac{x+y}{2} \right\| $$
defines a symmetric on $S_X$. It is elementary to observe that the $\rho$-diameter of a slice $\{x \in S_X: x^*(x)>1-\xi\}$ is bounded above by $\xi$ if $\|x^*\|=1$. Moreover, if the norm is $\sigma(X,F)$-LUR then $\rho$ 
semi-metrizes $(S_X,\sigma(X,F))$. Thus statement $(ii)$ is satisfied in both theorems.\\
If $(ii)$ holds, for every $n \in {\Bbb N}$ the family ${\mathcal F}_n$ of all the 
$\sigma(X,F)$-open slices of $S_X$ of $\rho$-diameter less than $1/n$ is a covering of 
$S_X$. Let $x,y \in S_X$ with $x \not = y$, then $y \not \in {\rm St}({x,\mathcal F}_n)$ whenever $n > 1/\rho(x,y)$. So $S_X$ has a ${\mathcal G}_\delta$-diagonal with $({\mathcal F}_n)_{n \in \N}$. Moreover, if $\rho$ semi-metrizes 
$(S_X,\sigma(X,F))$ then $({\mathcal F}_n)_{n \in \N}$ is a development. Indeed, trivially we have ${\rm St}({x,\mathcal F}_n) \subseteq B_\rho(x,1/n)$.\\
Now assume $(iii)$. Let $\| \cdot \|_k$ be the norm given by Theorem \ref{localisation} for the family of slices ${\mathcal F}_k$ and take $|\!|\!| x |\!|\!|^2=\sum_{k=1}^{+\infty} 2^{-k}\|x\|_k^2$. Assume that we are given points such that
$$ \lim_{n\rightarrow+\infty} (2|\!|\!|x_n|\!|\!|^2+2|\!|\!|x|\!|\!|^2-|\!|\!|x_n+x|\!|\!|^2) =0. $$
Then for every $k \in \N$ there is $N_k \in \N$ such that if $n>N_k$ then $x_n \in {\rm St}(x,{\mathcal F}_k)$. If $({\mathcal F}_n)_{n \in \N}$ witnesses the ${\mathcal G}_\delta$-diagonal property, taking $x_n=y$ we get that $y=x$ and so $|\!|\!| \cdot |\!|\!|$ is a rotund norm. Moreover, if $({\mathcal F}_n)_{n \in \N}$ is a development, then we get the $\sigma(X,F)$-convergence of $x_n$ to $x$, therefore the norm would be $\sigma(X,F)$-LUR.\findemo

\section{Dual Banach spaces}

The possibility of changing slices by general open sets in dual Banach spaces endowed with the weak$^*$-topology is based on a result from \cite{Raj07} that can be easily adapted to covering language. We shall need some terminology.
Let $X$ be a locally convex space and let ${\mathcal U}$ be a family of open sets. 
We write $A \preceq {\mathcal U}$ if there exists $U \in {\mathcal U}$ such that $A \subseteq U$.
Denote  by ${\mathcal H}$ the family of open half-spaces of $X$.
Given $A \subseteq X$  take
\[ [ A ]'_{\mathcal U} := \{ x \in A: \forall H \in {\mathcal H}, x \in H \Rightarrow A \cap H \not \preceq {\mathcal U} \}.  \]
Note that $[ A ]'_{\mathcal U}$ is convex if $A$ is so. These operations can be iterated, indeed we set $[A]^{k+1}_{\mathcal U} := [[A]^{k}_{\mathcal U}]'_{\mathcal U}$ and 
$[A]^{\omega}_{\mathcal U} := \bigcap_{n=1}^{\infty} [A]^{n}_{\mathcal U}$. Finally, we recall that ${\rm ext}(A)$ denotes the set of extreme points of a convex set $A \subseteq X$.

\begin{prop}\label{extremo}
If $A \subseteq X$ is compact and convex and let ${\mathcal U}$ be a family of open sets, then ${\rm ext}([A]^{\omega}_{\mathcal U}) \subseteq A \setminus \bigcup {\mathcal U}$.
\end{prop}

\noindent
{\bf Proof.} Suppose it is not the case. If $x \in {\rm ext}([C]_{\varepsilon}^{\omega})
\cap \bigcup {\mathcal U}$, there is $U \in {\mathcal U}$ such that $x \in U$. As $x$ is extreme, by Choquet's Lemma, there is $H \in {\mathcal H}$ such that $x \in H$ and  $\overline{H} \cap [A]^{\omega}_{\mathcal U} \subseteq U $.
Since $[A]^{\omega}_{\mathcal U} = \bigcap_{n=1}^{\infty} [A]_{\mathcal U}^{n}$, by compactness
there is $n \in {\Bbb N}$ such that $\overline{H} \cap [A]_{\mathcal U}^{n} \subseteq U$. Then we have
$H \cap [A]_{\mathcal U}^{n} \preceq {\mathcal U}$, thus $x \not \in [A]_{\mathcal U}^{n+1}$ which is a contradiction.\findemo

\begin{rema}
Note that if $B$ is a closed and convex set such that ${rm ext}([A]^{\omega}_{\mathcal U}) \cap B = \emptyset $, the half-space $H$ in the proof, and so all the ones that participates in the definition of $[ A ]'_{\mathcal U}$, can be taken such that $H \cap B = \emptyset$.
\end{rema}

\noindent
{\bf Proof of Theorem \ref{moore-dual}.} We already know, see 
\cite[Theorem~1.6]{FOR}, the equivalence of $(i)$, $(iv)$ and $(v)$. We have 
$(iv)$ implies $(iii)$ trivially, and $(iii)$ implies $(ii)$ just arguing like in the similar statements from Theorem \ref{moore-slice}. It will be enough to show that $(ii)$ implies $(i)$. Let $({\mathcal U}_k)_{k \in \N}$ be a development of 
$(S_X,w^*)$. We may assume, without loss of generality, that all the $w^*$-sets of $\bigcup_{k\in\N}\mathcal{U}_k$ skip $(1/2)B_X$. 
Fix $k \in \N$ and consider the sets $C_{k,j}=[B_X]_{{\mathcal U}_k}^{j+1}$ and $C_{k,1}=B_X$. By construction all those sets are $w^*$-closed convex symmetric and contains $(1/2)B_X$. 
Note as well that ${\rm ext}([B_X]^{\omega}_{{\mathcal U}_k}) \subseteq B_X \setminus \bigcup {\mathcal U}_k$ and so $\bigcap_{j=1}^\infty C_{k,j} \cap S_X = \emptyset$.
For every $j \in \N$ consider the family of slices of $C_{k,j}$ given by
$$ {\mathcal F}_{k,j} = \{ C_{k,j} \cap H: H \in {\mathcal H}, C_{k,j} \cap H
\preceq {\mathcal U}_k \}. $$
Let $\| \cdot \|_{k,j}$ the norm given by Theorem \ref{localisation} and take
$$ |\!|\!| x |\!|\!|^2 = \sum_{k,j=1}^{+\infty} 2^{-k-j} \| x\|^2_{k,j}. $$
We claim that this norm is $w^*$-LUR. Indeed, given points such that
$$ \lim_{n\rightarrow+\infty} (2|\!|\!|x_n|\!|\!|^2+2|\!|\!|x|\!|\!|^2-|\!|\!|x_n+x|\!|\!|^2) =0, $$
then a similar formula holds for each norm $\| \cdot \|_{k,j}$, in particular for the original norm too. Thus $\lim_n \|x_n\|=\|x\|$. If $\|x\|=0$ the sequence is norm converging and there is nothing to prove. Otherwise, we may divide $x$, $x_n$ by their norms and so we may assume that they are on the unit sphere of the original norm, say $S_X$. Fix $V$ a $w^*$-open neighbourhood of $x$, then fix $k \in \N$ such that ${\rm St}(x, {\mathcal U}_k) \subseteq V$. Now fix $j \in \N$ such that $x \in C_{k,j} \setminus C_{k,j+1}$. That implies ${\rm St}(x, {\mathcal F}_{k,j}) \subseteq V$. Clearly $\|x\|_{k,j}=1$. Since $\lim_{n\rightarrow+\infty} \|x_n\|_{k,j}=1$ we may do a further change taking $x'_n = x_n/ \|x_n\|_{k,j}$. Obviously 
$$ \lim_{n\rightarrow+\infty} (2 \|x'_n\|_{k,j}^2+2\|x\|_{k,j}^2-\|x'_n+x\|_{k,j}^2) =0 $$
and $\lim_{n\rightarrow+\infty} \|x'_n - x_n \| =0$. Thus for $n$ large enough we have $x'_n \in V$ and so $w^*$-$\lim_{n\rightarrow+\infty} x_n= x$, which finishes the proof of the theorem.\findemo\\

\noindent
{\bf Proof of Theorem \ref{strict-dual}.} $(i) \Leftrightarrow (iii)$ was established in \cite[Theorem 2.8]{FOR}, however here we only need the easier implication which makes use of the same symmetric and arguments that the proof of Theorem \ref{moore-slice}. $(i) \Leftrightarrow (v)$ is proved in \cite{OST}.
Clearly we have $(iii) \Rightarrow (ii) \Rightarrow (iv)$. It remains to show that $(iv) \Rightarrow (i)$. The proof will use the construction of the proof of Theorem \ref{moore-dual}, that is, the sets $C_{k,j}$, the families of slices 
${\mathcal F}_{k,j}$ and the norms $\| \cdot \|_{k,j}$. The goal now is to show the rotundity of $ |\!|\!| x |\!|\!|^2 = \sum_{k,j} 2^{-k-j} \| x\|^2_{k,j}$.
Indeed, if 
$$ 2|\!|\!|x|\!|\!|^2+2|\!|\!|y|\!|\!|^2-|\!|\!|x+y|\!|\!|^2 =0 $$
then $\| x \|_{k,j}=\| y \|_{k,j}=\| x+y \|_{k,j}/2$ for all $k,j \in \N$. If $x \not = y$, there would exist $k \in \N$ such that ${\mathcal U}_k$ separates $\{x,y\}$. As all the norms take  the same value on $x$ and $y$, there is $j \in \N$ such that $\{x,y\} \subseteq C_{k,j} \setminus C_{k,j+1}$. But this would imply that $y \not \in {\rm St}(x,{\mathcal F}_{k,j})$ and so $\| x \|_{k,j}=\| y \|_{k,j}=\| x+y \|_{k,j}/2$ leads to a contradiction.\findemo

\end{document}